\documentclass[11pt]{amsart}
\usepackage{amsmath}
\usepackage{amssymb}
\usepackage{epsfig}

\theoremstyle{plain}
\newtheorem{thm}{Theorem}[section]
\newtheorem{lem}[thm]{Lemma}
\newtheorem{cor}[thm]{Corollary}
\newtheorem{prop}[thm]{Proposition}

\theoremstyle{definition}

\def \R {\mathbf{R}}
\def \Z {\mathbf{Z}}
\def \Sig{\Sigma}

\def \CM {\mathcal M}

\def \a {\alpha}
\def \b {\beta}
\def \g {\gamma}

\def \lam {\lambda}
\def \G {\Gamma}

\def \t {\tau}

\def \bd {\partial}
\def \x {\times}

\def \ssw {\text{SW}}
\def \sw {\mathcal{SW}}

\def \DD {\Delta}

\def \ssw {\text{SW}}
\def \sw {\mathcal{SW}}

\def \DD {\Delta}

\begin{document}

\baselineskip.525cm

\title{Nonsymplectic $4$-Manifolds with One Basic Class}
\author[Ronald Fintushel]{Ronald Fintushel}
\address{Department of Mathematics, Michigan State University \newline
\hspace*{.375in}East Lansing, Michigan 48824}
\email{\rm{ronfint@math.msu.edu}}
\thanks{The first author was partially supported NSF Grant DMS9704927 and
the second author by NSF Grant DMS9626330}
\author[Ronald J. Stern]{Ronald J. Stern}
\address{Department of Mathematics, University of California \newline
\hspace*{.375in}Irvine,  California 92697}
\email{\rm{rstern@math.uci.edu}}
\date {July 27, 1999}

\maketitle

\section{Introduction\label{Intro}}

In the past several years Zoltan Szabo \cite{S1,S2} and the authors
\cite{KL4M} have produced examples of simply connected irreducible
4-manifolds which do not admit a
symplectic structure with either orientation. We shall call such manifolds
{\it{nonsymplectic}}. Due to the nature of their construction, these
manifolds have
many basic classes. It is the purpose of this paper to construct families
of examples of
nonsymplectic 4-manifolds which (up to sign) have just one basic class.

The key to detecting that the manifolds of \cite{S1,S2} and of \cite{KL4M}
are not symplectic lies in the theorem of C. Taubes which states that the
Seiberg-Witten
invariant associated to the canonical class of a symplectic $4$-manifold is
$\pm1$.
Recall that the Seiberg-Witten invariant  $\ssw_X$ of a smooth closed
oriented $4$-manifold
$X$ with $b ^+>1$ is an integer valued function which is defined on the set of
$spin^{\, c}$ structures over $X$ (cf. \cite{W}). In case $H_1(X;\Z)$ has no
2-torsion, there is a natural identification of the
$spin^{\, c}$ structures of
$X$ with the characteristic elements of $H^2(X;\Z)$. In this case we view the
Seiberg-Witten invariant as
\[ \ssw_X: \lbrace k\in H^2(X,\Z)|k\equiv w_2(TX)\pmod2)\rbrace
\rightarrow \Z. \] The Seiberg-Witten invariant $\ssw_X$ is a smooth
invariant whose sign depends on an orientation of
$H^0(X;\R)\otimes\det H_+^2(X;\R)\otimes \det H^1(X;\R)$. If $\ssw_X(\b)\neq
0$, then $\b$ is called a {\it{basic class}} of $X$. It is a fundamental
fact that the
set of basic classes is finite.  If $\b$ is a basic class, then so is $-\b$
with
\[\ssw_X(-\b)=(-1)^{(\text{e}+\text{sign})(X)/4}\,\ssw_X(\b)\] where
$\text{e}(X)$ is the Euler number and $\text{sign}(X)$ is the signature of $X$.
Because of this, we shall say that $X$ {\it{has $n$ basic classes}} if the
set $\{\b\,
|\,\ssw_X(\b)\ne0\}/ \{\pm1\}$ consists of $n$ elements.

There are abundant examples of $4$-manifolds with one basic class. Minimal
nonsingular
algebraic surfaces of general type have one basic class (the canonical class)
\cite{W}. The authors have constructed many examples of minimal symplectic
manifolds
with one basic class and $\chi -3\le {c_1}^2 < 2\chi -6$, where $\chi =
(b^++1)/2$.
(These manifolds cannot admit complex structures due to the geography of
complex
surfaces.) However, the examples produced below are the first nonsymplectic
manifolds
with one basic class.

As in \cite{KL4M} we need to view the Seiberg-Witten invariant as a Laurent
polynomial.
To do this, let
$\{\pm \b_1,\dots,\pm \b_n\}$ be the set of nonzero basic classes for $X$.
We may then view the Seiberg-Witten invariant of $X$ as the `symmetric'
Laurent polynomial
\[\sw_X = b_0+\sum_{j=1}^n
b_j(t_j+(-1)^{(\text{e}+\text{sign})(X)/4}\,t_j^{-1})\]
where  $b_0=\ssw_X(0)$, $b_j=\ssw_X(\b_j)$ and $t_j=\exp(\b_j)$.
The examples of \cite{S1,S2} and of \cite{KL4M} are obtained by producing
$4$-manifolds
whose Seiberg-Witten Laurent polynomial $\sw_X$ has as a factor a nonmonic
(symmetrized)
Alexander polynomial of a knot or link. Taubes' result is then used to show
that $X$
cannot have a symplectic structure. It is not difficult to see that any
nonsymplectic
manifold (with $b^+>1$) which can be constructed by the techniques of
\cite{KL4M} (as
explained in \cite{KL4M}, this includes the examples of Szabo) must have
more than one
basic class.

Whereas the examples of \cite{KL4M} are constructed by surgeries on
embedded tori of
self-intersection 0, the examples presented here arise from surgeries on
higher genus
surfaces. These examples are described in the next section.

\medskip

\section{A new family of 4-manifolds\label{family}}

We begin by recalling the construction of \cite{KL4M}. Suppose that we are
given a smooth
simply connected oriented $4$-manifold $X$ containing an essential smoothly
embedded torus
$T$ of self-intersection 0. Suppose further that $\pi_1(X\setminus T)=1$
and that $T$ is
contained in a cusp neighborhood. Let $K\subset S^3$ be a smooth knot and
$M_K$ the
3-manifold obtained from 0-framed surgery on $K$. The meridional loop $m$
to $K$ defines a
1-dimensional homology class $[m]$ both in $S^3\setminus K$ and in $M_K$.
Denote by $T_m$
the torus $S^1\x m\subset S^1\x M_K$. Then $X_K$ is defined to be the fiber
sum
\[ X_K= X\#_{T=T_m} S^1\x M_K = (X\setminus N(T))\cup(S^1\x(S^3\setminus
N(K)), \]
where $N(T)\cong D^2\x T^2$ is a tubular neighborhood of $T$ in $X$ and
$N(K)$ is a
neighborhood of $K$ in $S^3$. If $\lam$ denotes the longitude of $K$
($\lam$ bounds a
surface in $S^3\setminus K$) then the gluing of this fiber sum identifies
$\{ {\text{pt}}
\}\x \lam$ with a normal circle to $T$ in $X$. The main theorem of
\cite{KL4M} asserts
that $X_K$ is homeomorphic to $X$, and
\[ \sw_{X_K}=\sw_X\cdot \DD_K(t)\]
where $\DD_K$ is the symmetrized Alexander polynomial of $K$ and
$t=\exp(2[T])$.

To begin our construction, take $X$ to be the $K3$-surface (which has
$\sw_X=1$) and let
the torus $T$ be a smooth fiber of an elliptic fibration on $X$. The pair
$(X,T)$
satisfies the hypotheses of  \cite{KL4M}; so for any knot $K$ we get a homotopy
$K3$-surface $X_K$ whose Seiberg-Witten invariant is $\sw_{X_K}=\DD_K(t)$. The
$K3$-surface,
$X$, has a section (to the elliptic fibration) which is a smoothly embedded
2-sphere $S$
of self-intersection
$-2$, and $[S]\cdot [T]=1$. The sum $[S]+[T]$ is represented by a smooth
torus $\Sig$ with
$[\Sig]^2=0$ and $[\Sig]\cdot[T]=1$.

Suppose that the knot $K$ has genus $g$. In the construction of
$X_K$ we have replaced a 2-disk in $S$ (normal to $T$) with a punctured
surface of genus
$g$. Thus $X_K$ contains a genus $g$ surface $S'$ satisfying
$[S']\cdot[S']=-2$ and
$[S']\cdot [T] = 1$. Consider another smooth fiber $T'$ of the elliptic
fibration of
$(X\setminus N(T))\subset X_K$. Then $T'+S'$ is a singular curve with one
double point,
which can be smoothed to give an embedded surface $\Sig'$ of genus $g+1$
representing the
homology class $[\Sig']=[T']+[S']$. Thus $[\Sig']^2=0$ and $[\Sig']\cdot
[T]=1$.

Next, let $K'$ denote the left-handed trefoil knot in $S^3$. Since $K'$ is
a fibered genus 1
knot, the 4-manifold $S^1\x M_{K'}$ is a smooth $T^2$-fiber bundle over
$T^2$. Forming the
fiber sum of $g+1$ copies of $S^1\x M_{K'}$, we obtain
\[ \begin{array}{ccc}
F=T^2 \longrightarrow & Y&\negthickspace =\ S^1\x M_{K'}\#_F\cdots\#_F
S^1\x M_{K'}\\
&\Big\downarrow&\\
& C_0&
\end{array}\]
where $C_0$ is a genus $g+1$ surface. Furthermore, $S^1\x M_{K'}$ is a
symplectic manifold
\cite{Th}.  Notice that there is a section $C\subset Y$ of
the fibration given by the connected sum of the $g+1$ tori $T_{m_i}$.

Generally, suppose that we are given symplectic 4-manifolds $A$ and $B$ and
that
$A\#_NB$ is their symplectic fiber sum along a symplectic torus of
self-intersection
$0$. The adjunction formula implies that the canonical class $K_A$ is
orthogonal to
$[N]$, as is $K_B$. The canonical class of $A\#_NB$ is then
$K_{A\#_NB}=K_A+K_B+2[N]$
({cf.\cite{Gompf}). Apply this fact to $X_{K'}=K3\#_{T=T_m}S^1\x M_{K'}$. Since
$\sw_{X_{K'}}=\DD_{K'}(t)=t-1+t^{-1}$ (where $t=\exp(2[T_m]$), we see that
$ K_{X_{K'}}= 2T_m = K_{K3}+K_{S^1\x M_{K'}}+2[T_m]=K_{S^1\x M_{K'}}+2[T_m]$.
Hence $K_{S^1\x M_{K'}}=0$ and so $c_1(K_{S^1\x M_{K'}})=0$. Now apply the
fact $g$ more
times, this time fiber-summing along $F$. It follows that $Y$ is a
symplectic 4-manifold
with $c_1(Y)=-2g[F]$.
(Here, we identify $[F]$ with its Poincar\'e dual.)

Our example, corresponding to the genus $g$ knot $K$ is
\[ Z_K=X_K\#_{\Sig'=C}Y.\]
We perform this fiber sum so that $Z_K$ is a spin 4-manifold \cite{Gompf}.

\begin{prop} The manifold $Z_K$ is simply connected. \end{prop}
\begin{proof} The fundamental group $\pi_1(X_K\setminus\Sig')$ is normally
generated by a
boundary circle of a normal disk to $\Sig'$. Since $[\Sig']\cdot[T]=1$, we
may assume that
this circle lies on a copy $\{\text{pt}\}\x T$ in the boundary $\bd D^2\x
T=\bd N(T)$.
We claim that there are generators $\lam_1$, $\lam_2$, for $\pi_1(T)$ which
bound
vanishing cycles (disks of self-intersection $-1$) in $X\setminus (S\cup
T')$. (Note that
here we are identifying $X\setminus T$ with $X_K\setminus T$.) This claim
can be seen to
be true inside a $K3$ nucleus, i.e. in a regular neighborhood of the union
of $S$ and a
cusp fiber. A Kirby calculus diagram for the nucleus is given in the
Figure~1 below.

The section $S$ is the union of the disk spanned by the circle
labelled `$-2$' and the core of the 2-handle which is attached to it.
The torus $T$ is obtained as follows. The circle labelled `0' bounds a disk
$D$ which is
punctured in two points by one of the dotted circles. Remove a pair of
disks from $D$ at
these intersection points, and connect the boundaries of these disks with
an annulus
which surrounds the path $\g$ in the diagram. The torus $T$ is the union of
the
twice-punctured $D$ and this annulus. We can see that the loop $\a$ of the
diagram lies
on $T$, and it is easy to deform $\b$ to also lie on $T$. (When this is
done, $\a$ and
$\b$ will intersect in a point.)
Thus $H_1(T)$ is generated by the classes represented by the loops $\a$ and
$\b$. The vanishing cycles are the cores of the  ($-1$)-framed 2-handles which
are attached to $\a$ and $\b$. This proves the claim.

\centerline{\unitlength 1cm
\begin{picture}(4,6.5)
\put (1,5){\oval(1,1)[l]}
\put (1.2,4.7){\line(0,1){1.15}}
\put (2.9,4.7){\line(0,1){1.15}}
\put (3,5){\oval(1,1)[r]}
\put (1,4.485){\line(1,0){2}}
\put (.71,5.85){\oval(1,1)[t]}
\put (3.41,5.85){\oval(1,1)[t]}
\put (0.2,5.85){\line(0,-1){3.85}}
\put (3.9,5.85){\line(0,-1){4}}
\put (1.2,4.25){\line(0,-1){1.2}}
\put (2.9,4.25){\line(0,-1){.5}}
\put (2.06,2){\oval(3.7,1)[bl]}
\put (1.55,3.035){\oval(.675,.55)[bl]}
\put (1.5,2.75){\line(1,0){.65}}
\put (1.96,2.9){\oval(.5,.75)[t]}
\put (2.2,2.9){\line(0,-1){1.75}}
\put (1.96,1.25){\oval(.5,.75)[b]}
\put (1.7,2.65){\line(0,-1){.6}}
\put (1.7,1.85){\line(0,-1){.3}}
\put (1.7,1.4){\line(0,-1){.2}}
\put (1.7,2.2){\oval(.5,.5)[b]}
\put (1.6,2.2){\oval(.3,.3)[tl]}
\put (1.8,2.2){\oval(.3,.3)[tr]}
\put (2.9,3.85){\oval(.5,.5)[b]}
\put (2.8,3.85){\oval(.3,.3)[tl]}
\put (3,3.85){\oval(.3,.3)[tr]}
\put (3.2,3.75){\footnotesize$-1$}
\put (3.05,3.45){\tiny$\a$}
\put (1,2.1){\footnotesize$-1$}
\put (1.35,1.77){\tiny$\b$}
\put (2.9,3.5){\line(0,-1){.5}}
\put (2.25,2.75){\line(1,0){.35}}
\put (2.57,3.035){\oval(.675,.55)[br]}
\put (2.25,1.87){\oval(3.32,.75)[br]}
\put (2.1,2.1){$\bullet$}
\put (3.8,3.75){$\bullet$}
\put (2,4.2){\footnotesize$0$}
\put (2,5.5){\oval(.5,.75)[r]}
\put (1,5.49){\line(1,0){.1}}
\put (1.3,5.49){\line(1,0){.85}}
\put (2.3,5.49){\line(1,0){.5}}
\put (2,5.625){\oval(.35,.5)[tl]}
\put (2,5.375){\oval(.35,.5)[bl]}
\put (2,5.9){\footnotesize$-2$}
\put (0.95,3.5){\tiny$\g$}
\put (1.5,.25){Figure 1}
\end{picture}}

This means that  $\pi_1(\{\text{pt}\}\x T)\to\pi_1(X\setminus\Sig)$ is the
zero map;
hence $\pi_1(\{\text{pt}\}\x T)\to\pi_1(X_K\setminus\Sig')$ is also the
zero map.
However, $\pi_1(X_K\setminus\Sig')$  is normally generated by the image of
this map; so
$X_K\setminus\Sig'$ is simply connected. Thus
\[ \pi_1(Z_K)=\pi_1(X_K\setminus\Sig')*_{\pi_1(C\x S^1)}\pi_1(Y\setminus C)
  = \pi_1(Y\setminus C)/\pi_1(C\x S^1). \]

Because $\pi_1(S^3\setminus K')$ is normally generated by the meridian $m$,
$\pi_1(S^1\x (S^3\setminus K'))=\pi_1((S^1\x M_{K'})F)$ is normally
generated by
$\pi_1(S^1\x m)$. An inductive application of Van Kampen's theorem shows
that $\pi_1(Y)$
is normally generated by $\pi_1(S^1\x m \#\cdots\# S^1\x m)=\pi_1(C)$. Thus
$\pi_1(Y\setminus C)$ is normally generated by $\pi_1(C\x S^1)$, and so
$\pi_1(Z_K)=1$.
\end{proof}

\medskip
\section{The Seiberg-Witten invariants of $Z_K$\label{SW}}

Consider first $H_2(Z_K)$. There is an important class $\t\in H_2(Z_K)$
constructed as
follows. In the construction of $Z_K$, the boundary of a tubular
neighborhood $N(\Sig')$
of $\Sig'$ in $X_K$ is identified with the boundary of a tubular
neighborhood $N(C)$
of $C$ in $Y$. Fix a fiber $F$ of $Y$, and let $F_0=F\setminus(F\cap
\text{int} N(C))$.
There is torus $T''$ which is a smooth fiber of the elliptic fibration of
$X\setminus\{ T\cup T'\}\subset X_K$ and such that if
$T''_0=T''\setminus(T''\cap \text{int} N(\Sig'))$, then $\bd T''_0=\bd F_0$
in $Z_K$. Let
$\t$ denote the class $\t= [T''_0\cup F_0]$. Note that $\t$ is represented
by an
embedded surface of genus 2, $\t^2=0$, and $\t\cdot[\Sig']=1$. Then
$H_2(Z_K)$ is generated
by the image of $H_2(Y\setminus C)$, of $H_2(X_K\setminus\Sig')$, and the
class $\t$. The
only other classes which could contribute to $H_2(Z_K)$ are the classes of
rim tori, i.e.
tori lying on $\bd N(\Sig')=\bd N(C)$ in $Z_K$ which have the form
$\xi\x\bd D^2$ where $D^2$ is a normal disk to $\Sig'$ (or to $C$). The
next lemma shows
that in fact they are all trivial.

\begin{lem}\label{rim} Each rim torus is homologically trivial in $Z_K$.
\end{lem}
\begin{proof} A rim torus on $\bd N(C)$ has the form $\xi\x\bd D^2$, for
some loop $\xi$
on $C$. Recall that there is a fiber bundle $\varphi: Y\to C_0$ with fiber
$F$. Let
$Q=\varphi^{-1}(\xi)\subset Y$. We see that $\xi\x\bd D^2$ bounds the
3-chain $Q\setminus
(\xi\x\,\text{int}\,D^2)$ in $Y\setminus N(C)\subset Z_K$.\end{proof}

Before we prove our main theorem, we recall that a 4-manifold $W$ is said
to have
{\it{simple type}} if $\ssw_W(k)\ne 0$ implies that
\[ \dim\CM_W(k)= \frac14 (k^2-(3\,\text{sign}+2\,\text{e})(W))=0 \]
where $\CM_W(k)$ is the Seiberg-Witten moduli space. Write the symmetrized
Alexander polynomial of $K$ as
$\DD_K(t)=a_0+\sum\limits_{n=1}^da_n(t^n+t^{-n})$. We
call $d$ the {\it{degree}} of $\DD_K(t)$. We are assuming that the genus of
$K$ is $g$; so
$d\le g$. If $K$ is an alternating knot, for example, then $d=g$. Let us
say that the
Alexander polynomial of $K$ has {\it{maximal degree}} if $d=g$.

\begin{thm} Let $K$ be a knot in $S^3$ whose Alexander polynomial has
maximal degree.
Then $Z_K$ is of simple type and has (up to sign) a single basic class,
$k=2g\,\t+2[\Sig']$. Furthermore, $|\ssw_{Z_K}(k)|=a_d$, the top coefficient of
$\DD_K(t)$.
\end{thm}
\begin{proof} Let $U$ denote a nucleus in $X=K3$ which contains the fiber
$T$ and section
$S$ from the construction of $X_K$. We see that $(X\setminus U)\subset
(X_K\setminus\Sig')$. The homology $H_2(X\setminus U)\cong 2\,E_8\oplus
2\,H$, where the
negative definite $E_8$ forms are generated by the classes of embedded
spheres of
self-intersection $-2$, and the two hyperbolic pairs $H$ are each generated
by a torus
$T_i$ of self-intersection 0, and a sphere $S_i$ of self-intersection -2
which meet
transversely in a single point. The homology $H_2(X_K\setminus\Sig')$ is
generated by the
image of $H_2(X\setminus U)$ together with $[\Sig']$ and the classes of rim
tori.

Next consider $Y = S^1\x M_{K'}\#_F\cdots\#_F S^1\x M_{K'}$ ($g+1$ copies)
where $F$ is
the torus fiber of the fibration of $S^1\x M_{K'}$ over the torus. Each
$S^1\x M_{K'}$
has the homology of $S^2\x T^2$. Each fiber sum in the construction of $Y$
increases the
first betti number
$b_1$ by $2$ --- the base of the fibration has genus increased by $1$ ---
and $H_2$ is
increased by the addition of two hyperbolic pairs as follows: Choose a
standard basis
$\{ x_1, x_2\}$ for $H_1(F;\Z)$. For example, $x_1$ is represented by a
loop as shown in
Figure~2.  Each of the curves $x_i$ bounds a punctured torus $\G_i$ in
$M_K$. In
Figure~3, $x_1$ is isotopic to the pictured curve $x_1'$, and the punctured
torus is composed
of the twice-punctured disk which $x_1'$ bounds, together with a 1-handle
which pipes around
the knot. Let $x_i''$ be a push-off of $x_i$ in $F$. Then the linking
number of $x_i$ and
$x_i''$ is $+1$. (Here we are using the fact that $K'$ is a left-hand
trefoil knot.) This means
that the self-intersection number of $\G_i$ (as a surface in $M_K\x I$,
say), keeping its
boundary in
$F$, is $+1$.

\centerline{\unitlength 1cm
\begin{picture}(12,6)
\put (2,3){\oval(3,4)[l]}
\put (3,3){\oval(3,4)[r]}
\put (2,5){\line(3,-4){1}}
\put (3,3.66){\line(-3,-4){.4}}
\put (3,5){\line(-3,-4){.4}}
\put (2,3.66){\line(3,4){.4}}
\put (2,3.66){\line(3,-4){1}}
\put (2,2.33){\line(3,4){.4}}
\put (2,2.33){\line(3,-4){1}}
\put (3,2.33){\line(-3,-4){.4}}
\put (2,1){\line(3,4){.4}}
\put (2.4,3.66){\oval(1.5,1.4)[l]}
\put (2.6,3.66){\oval(1.5,1.4)[r]}
\put (1.5,4.25){\small{$x_1$}}
\put (4.65,4.4){\small{$0$}}
\put (1.85,.5){Figure 2}
\put (9,3){\oval(3,4)[bl]}
\put (7.49,3){\line(0,1){.125}}
\put (10,3){\oval(3,4)[r]}
\put (9,5){\line(-1,0){.85}}
\qbezier (7.49,4.35)(7.6,4.9)(8.15,5)
\put (9,5){\line(3,-4){1}}
\put (10,3.66){\line(-3,-4){.4}}
\put (10,5){\line(-3,-4){.4}}
\put (9,3.66){\line(3,4){.4}}
\put (9,3.66){\line(3,-4){.265}}%%%
\put (9,2.33){\line(3,4){.4}}
\put (9,2.33){\line(3,-4){1}}
\put (10,2.33){\line(-3,-4){.4}}
\put (10,2.33){\line(-3,4){.6}}%%
\put (9,1){\line(3,4){.4}}
\put (7.1,3.2){\line(1,0){2.3}}
\put (7.65,4){\line(1,0){1.5}}
\put (7.1,4){\line(1,0){.25}}
\qbezier (9.4,3.2)(10,3.6)(9.4,4)
\qbezier (7.1,3.2)(6.5,3.6)(7.1,4)
\put (7.49,4.35){\line(0,-1){1}}%%%%%%
\put (6.4,3.75){\small{$x_1'$}}
\put (11.65,4.4){\small{$0$}}
\put (8.85,.5){Figure 3}
%\put (11.75,5.35){\line(-1,0){.25}}
\end{picture}}

 Thus in
$S^1\x M_{K'}\#_F S^1\x M_{K'}$ one produces genus 2 surfaces $S_1'$,
$S_2'$, of
self-intersection $+2$ which are formed from pairs of these tori. Let $T_1'$,
$T_2'$ be the rim tori corresponding to $x_2$, $x_1$ (reversed on purpose).
Then in
$H_2(S^1\x M_{K'}\#_F S^1\x M_{K'};\Z)$ two
hyperbolic pairs are generated by the pairs $\{ [S_i'], [T_i']\}$. Each
further fiber
sum adds two such hyperbolic pairs to $H_2$. It follows that
$H_2(Y\setminus C)$, is generated
by the
$[S_i']$,
$[T_i']$, and the section class $[C]$.

Using our observations above, if $k$ is a basic class of $Z_K$ we can write
\[ k= a\t+b[\Sig']+\b+\sum_{i=1}^2 m_i[T_i]+n_i[S_i] + \sum_{j=1}^{2g} t_j[T_j']+s_j[S_j']\]  where $a>0$ and $\b\in 2\,E_8\subset H_2(X\setminus U)$.T
he adjunction inequality (see e.g. \cite{MST}) states that if $k$ is a
basic class and
$B$ is an embedded surface of genus $g_B$ and self-intersection $[B]^2\ge
0$ then
\[ 2g_B-2\ge [B]^2 + |k\cdot[B]|. \]
In particular, this implies that if the self-intersection of $B$ is
$[B]^2=2g_B-2$,
then any basic class
$k$ must be orthogonal to it: $k\cdot [B]=0$. Since $T_i$ is a torus of
self-interesection
$0$, it follows that
$n_i=k\cdot [T_i]=0$, and, since $[S_i]+[T_i]$ is also represented by an
embedded torus of
self-intersection
$0$, $m_i+n_i=k\cdot ([T_i]+[S_i])=0$. The same argument applies to show
that $s_j=0=t_j$
for each $j=1,\dots,2g$. Thus
$k= a\t+b[\Sig']+\b$. Apply the adjunction inequality to the genus $g+1$
surface $\Sig'$
and the genus 2 surface representing $\t$ to obtain:
\begin{equation}\tag{*} a=k\cdot [\Sig']\le 2g, \hspace{.5in} |b|=
|k\cdot\t|\le 2.
\label{*}\end{equation}

Because $k$ is a basic class, $\dim\CM_{Z_K}(k)\ge0$, hence
\[ 0\le k^2-(3\,\text{sign}+2\,\text{e})(Z_K). \]
Since $(3\,\text{sign}+2\,\text{e})(X_K)=0$, and
$(3\,\text{sign}+2\,\text{e})(Y)=0$, it is
easy to check that $(3\,\text{sign}+2\,\text{e})(Z_K)=8g$. Furthermore,
$\b$ lies in
the negative definite space $2\,E_8$; so if $\b\ne0$ then
\[ 0\le 2ab+\b^2-8g<2ab-8g\le 8g-8g =0.\]
This contradiction implies that $\b=0$; so $k= a\t+b[\Sig']$. Any of the
($-2$)-spheres
generating the $E_8$'s is orthogonal to $k$; hence orthogonal to each basic
class of
$Z_K$. It now follows from \cite{FS} that $Z_K$ has simple type;
$\dim\CM_{Z_K}(k)=0$. Thus we have
\[2ab=k^2=(3\,\text{sign}+2\,\text{e})(Z_K)=8g.\]
It now follows from \eqref{*} that $a=2g$ and $b=2$, as claimed.

Finally, we apply a theorem of Morgan, Szabo, and Taubes to calculate
$\ssw_{Z_K}(k)$.
Since $k\cdot\Sig'=2g$, \cite{MST} applies to give
\[\ssw_{Z_K}(k) =\sum\ssw_{Z_K}(k+2[R])=
\pm\ssw_{X_K}(2g\,T)\cdot\ssw_Y(2g\,F)\]
where the the sum is taken for all distinct classes $k+2[R]$ for $R$ a rim
torus.
Thus the first equality follows from Lemma \ref{rim} which shows that each
$[R]=0$ in
$H_2(Z_K;\Z)$.
Now $Y$ is a symplectic manifold with $c_1(Y)=-2g\,F$; so \cite{T1} implies
that
$\ssw_Y(2g\,F)=\pm1$. Since we are assuming that $g=d$, \cite{KL4M}
implies that $\ssw_{X_K}(2g\,T)=a_d$, and this completes the proof.
\end{proof}

We remark that in case the Alexander polynomial of the knot $K$ does not
have maximal
degree, the above proof shows that $\sw_{Z_K}=0$; this provides potential
examples of
simply connected irreducible 4-manifolds with trivial Seiberg-Witten
invariants.

\begin{cor} Let $K$ be a knot in $S^3$ whose Alexander polynomial is not
monic and has
maximal degree. Then $Z_K$ is a nonsymplectic 4-manifold with one basic
class.\end{cor}
\begin{proof} The hypothesis implies that the only nontrivial
Seiberg-Witten invariant of
$Z_K$ has value $\pm a_d\ne\pm1$; so \cite{T1} implies that $Z_K$ has no
symplectic
structure. Since $Z_K$ contains an embedded sphere of self-intersection
$-2$, neither does
$Z_K$ with its opposite orientation admit a symplectic structure.
\end{proof}

We close with a comment concerning the geography of our construction. If
the genus of
$K$ is $g$ then $c(Z_K)\equiv(3\,\text{sign}+2\,\text{e})(Z_K)=8g$ and
$\chi(Z_K)\equiv \frac{b^++1}{2}(Z_K)=g+2$; so all these manifolds lie on
the line of
slope $8$ passing through $(2,0)=(\chi(K3),c(K3))$ in the $(\chi,
c)$-plane. We could
similarly perform our construction starting with $X=E(2n)$, the minimal
elliptic surface
without multiple fibers and with $\chi=2n$. We then obtain the same result
for the
lattice points $(3n+g-1, 8(g+n-1))$, i.e. for the lattice points on the
line of slope $8$
through
$(2n,0)=(\chi(E(2n)),c(E(2n)))$.

\end{document}